
\input amstex
\magnification = \magstep1
\hsize = 6.25 truein
\vsize = 20 truecm
\baselineskip .22in
\define\Lap{\varDelta}
\define\del{\partial}
\define\a{\alpha}
\predefine\b{\barunder}
\redefine\b{\beta}

\predefine\d{\dotunder}
\redefine\d{\delta}
\define\th{\theta}
\define\e{\epsilon}

\predefine\o{\orsted}
\redefine\o{\omega}

\define\cin{\Cal C^{\infty}}
\define\RR{\Bbb R}
\define\CC{\Bbb C}
\define\cri{\frac{n+2}{n-2}}
\define\cl{\frac{4}{n-2}}
\define\co{\frac{n-2}{4(n-1)}}
\define\Le{L_{\e}}

\define\com{\frac{n+2}{4(n-1)}}
\define\La{\Lambda}
\define\ML{\Cal M_{\La}}
\define\Mk{\Bbb M_k}
\define\Ck{\Cal C_{k}}


\documentstyle{amsppt}
\NoRunningHeads
\topmatter
\title Connected sum constructions for constant scalar curvature
metrics \endtitle
\author Rafe Mazzeo${}^{(\dag)}$, Daniel Pollack${}^{(\ddag)}$ and
Karen Uhlenbeck${}^{(\star)}$ \endauthor
\affil Stanford University ($\dag$), University of Chicago ($\ddag$)
and University of Texas at Austin ($\star$) \endaffil
\dedicatory Dedicated to Louis Nirenberg on the occasion
of his 70$^{th}$ birthday\enddedicatory
\abstract
We give a general procedure for gluing together possibly
noncompact manifolds of constant scalar curvature which satisfy
an extra nondegeneracy hypothesis. Our aim is to provide a
simple paradigm for making `analytic' connected sums.
In particular, we can easily construct complete metrics of constant
positive scalar curvature on the complement of certain configurations
of an even number of points on the sphere, which is a special
case of Schoen's \cite{S1} well-known, difficult construction.
Applications of this construction produces metrics with
prescribed asymptotics. In particular, we produce metrics with
cylindrical ends, the simplest type of asymptotic behaviour.
Solutions on the complement of an infinite number of points
are also constructed by an iteration of our construction.
\endabstract
\thanks Research supported in part by ($\dag$) NSF Young Investigator
Award, the Sloan Foundation and NSF Grant \# 9303236, ($\ddag$) NSF
Postdoctoral Research Fellowship
and ($\star$) the Sid Richardson and O'Donnell Foundations. \endthanks

\address{Department of Mathematics, Stanford University, Stanford, CA 94305}
\endaddress
\email{mazzeo\@math.stanford.edu}
\endemail

\address{Department of Mathematics, University
of Chicago, Chicago, IL 60637}
\endaddress
\email{pollack\@math.uchicago.edu}
\endemail

\address{Department of Mathematics, University
of Texas, Austin, TX 78712}
\endaddress
\email{uhlen\@math.utexas.edu}
\endemail

\endtopmatter
\document

\specialhead I. Introduction \endspecialhead

It is now a well-entrenched procedure in geometric analysis
to construct new solutions to nonlinear PDE by gluing together
known solutions: an approximate solution is constructed, then
perturbed to an exact solution using analytic methods.
One of the early spectacular instances of this is Taubes'
patching of instantons \cite{T}. More recent
instances are too numerous to list here.

The geometric problem we wish to examine here is the possibility of
gluing together manifolds of constant scalar curvature satisfying
a certain nondegeneracy condition to obtain a new constant scalar
curvature metric on the connected sum. The precise notion of nondegeneracy
will be given in \S2 below, but when the manifolds are compact, possibly
with boundary, it coincides with the invertibility of the Jacobi operator
(which is weaker than stability). The main result is:

\proclaim{Theorem}
Let $(X_1,g_1)$ and $(X_2,g_2)$ be any two manifolds, possibly
with boundary, with complete metrics $g_1$, $g_2$ of constant scalar
curvature $n(n-1)$. Suppose also that the metrics $g_i$ satisfy
the nondegeneracy condition (2.12) and either (2.15) or (2.16-17) below.
Then for any points $p_i \in X_i$,
the connected sum $X_1 \#_\e X_2$ obtained by excising small $\e$-balls
around the $p_i$ and identifying boundaries, carries a complete
nondegenerate metric $g_\e$ of constant scalar curvature $n(n-1)$.
\endproclaim

The problem of gluing nondegenerate
compact constant scalar curvature manifolds has already been studied by
Joyce \cite{J}, so our primary interest is with noncompact manifolds.
Furthermore, we will focus exclusively on manifolds with constant
positive scalar curvature (CPSC). The simplest of these that we wish to
treat as `summands' to be glued are
the Delaunay metrics on the complement of two points in the
sphere $S^n$.  These are conformally equivalent to
elements of an explicit one-parameter
family of rotationally symmetric metrics interpolating between
the cylinder and an infinite bead of spheres strung out along
a common axis. These metrics satisfy the nondegeneracy
condition, and the metrics we construct on the connected
sum of any finite (or even infinite) number of these are conformally flat,
hence may be uniformized and regarded as complete metrics
on $S^n \setminus \La$, where $\Lambda$ is a discrete
collection of points of even cardinality. Complete CPSC metrics
on $S^n \setminus \La$, for $\La$ finite, were
originally constructed in Schoen's well-known and difficult
paper \cite{S1}. One motivation for our construction is
to provide a simple proof of a special case of his result.
The solutions of this type which we construct here are called
dipole metrics, because their singular sets $\La$ are
widely separated pairs of closely spaced points.

One contribution of this paper is a general formulation
of nondegeneracy of the scalar curvature operator on manifolds
with complete CPSC metrics. The importance of the nondegeneracy
hypothesis is clear, for example, in the analysis of the moduli space
$\ML$ of complete CPSC metrics on $S^n \setminus
\La$, where $\La$ is a submanifold.  When any
component of $\La$ has positive dimension, $\ML$ is infinite dimensional,
\cite{MPa}, but when $\Lambda$ is a finite set of points,
$\ML$  is a finite dimensional real analytic set, \cite{MPU}.
If $g \in \ML$ is a nondegenerate solution, then in a neighbourhood of
it, $\ML$ is a smooth (in fact, real analytic) manifold.
In particular, when $\La$ is finite,
this neighbourhood is of dimension equal to the cardinality of $\La$.
Unfortunately, we had previously been unable to establish the
nondegeneracy of the solutions constructed in \cite{S1}, so a number
of simple statements about the moduli space theory remained
hypothetical.  The nondegeneracy of our dipole solutions clarify many of
these moduli space issues.  This is discussed further in $\S 4$.
Even when $\La$ is positive dimensional, nondegeneracy
has important ramifications for the Dirichlet problem
parametrizing the infinite dimensional moduli space $\ML$.

It is possible, in certain circumstances, to glue metrics
not satisfying the nondegeneracy conditions. The main
instance is Schoen's construction \cite{S1}
mentioned earlier,  cf. also \cite{P1} for an analogous
construction in the compact case. The summands
in these constructions are standard spheres, for which
the Jacobi operator is definitely not invertible.
In a forthcoming paper \cite{MPa2}, a new and simpler
proof of Schoen's theorem will be given; the simplification
relies on the observation that here too there is an
underlying nondegenerate gluing procedure.

It is well-known that there are close relationships between
the problems concerning CPSC metrics and constant mean
curvature (CMC) surfaces in $\RR^3$. Indeed, closely related
in form to \cite{S1}, but substantially different in many
technical details, is Kapouleas' famous construction \cite{K} of
CMC surfaces, both compact and noncompact, in $\RR^3$.
It is possible to adapt the ideas here to construct
noncompact CMC surfaces; these surfaces are topologically
identical, but geometrically quite different from many of the
surfaces obtained by Kapouleas.  Because of the simplicity
of the CPSC construction, relative to that for CMC surfaces,
we defer the CMC construction to a subsequent paper.

In \S2 we first discuss the main examples of nondegenerate CPSC metrics
and then, motivated by these examples, give an abstract definition of
nondegeneracy. Next, in \S3, we use this to prove the main gluing theorem.
In \S4 we apply this to the special case where the summands in the gluing
construction are Delaunay metrics on the cylinder, and the
ramifications of this theorem for the moduli space theory
of \cite{MPU}. We also introduce here the `unmarked moduli space'
of CPSC metrics on the complement of any collection of
$k$ distinct points in $S^n$ and prove, analogously to \cite{MPU},
that it is a real analytic set; finally, we relate the nondegeneracy
of elements in this unmarked moduli space to
the problem of showing that solutions other than the ones obtained in \S3
by the gluing theorem are nondegenerate.

\specialhead II. Nondegeneracy: examples and definitions \endspecialhead

In this section we set up the notation used for the
remainder of the paper and then give a precise definition of the
nondegeneracy of solutions. We first motivate this definition
by describing in some detail the key examples which led to it.

Let $(M,g_0)$ be a fixed complete Riemannian manifold, which
we do not assume to have constant positive scalar curvature.
Suppose $g$ is a complete CPSC manifold conformal to $g_0$.
We express the conformal factor by writing $g = u^{\cl}g_0$.
Letting $R(g_0)$ and $R(g)$ denote the scalar curvatures of
$g_0$ and $g$ then it is well known that
$$
\Lap_{g_0}u - \co R(g_0) u + \co R(g) u^{\cri} = 0. \tag 2.1
$$
Denote the left side of this equation by $N_{g_0}(u)$. Much of
the analysis of CPSC metrics near $g$ revolves around the
linearization $L$ of $N_{g_0}$ at $u$:
$$
\gathered
L_{g_0} v = \left. \frac{\del\,}{\del t} \right|_{t=0} N_{g_0}(u + tv) \\
= \Lap_{g_0}v - \co R(g_0) v + \com R(g) u^{\cl} v.
\endgathered
\tag 2.2
$$
In a special case, where $g = g_0$ and $R(g) = n(n-1)$, this
operator takes the form
$$
L_g v = \Lap_g v + nv. \tag 2.3
$$
For convenience we let $L$ denote the linearization; whether it is
relative to $g_0$ or $g$ will be clear from the context.  Our interest
in this section is in the mapping properties of $L$.

When $M$ is compact, $L$ is self-adjoint, and in this case
it is said to be nondegenerate provided $0 \notin \text{spec\,}(L)$.
This is equivalent to either the injectivity or surjectivity of $L:
H^{s+2}(M) \rightarrow H^s(M)$ for any $s$.
Although it is the surjectivity that is used in the nonlinear
analysis, it is usually easier to check injectivity. For example,
it is clear that the sphere $S^n$ with its standard metric is degenerate
because $L$ annihilates the restrictions of linear functions on
$\RR^{n+1}$ to $S^n$.

When $M$ is noncompact, the precise formulation of nondegeneracy
is more subtle since in all the known examples $0$ {\it is}
in the spectrum of
$L$. Rather than exclude these, we must examine the mapping properties of
$L$ more closely.  Before stating the correct abstract formulation of
nondegeneracy, we present the two key examples motivating this definition.

\head Delaunay metrics \endhead

The punctured sphere $M = S^n \backslash \{p_1, p_2\}$ with its
standard metric has CPSC but is incomplete. However, it is conformal
to the complete CPSC product metric $g = dt^2 + d\th^2$ on
$\RR_t \times S^{n-1}_{\th}$.
There is a one-parameter family of complete CPSC metrics $g_{\e}$,
$0 < \e \le \bar u $ with  $\bar u = \left((n-2)/n\right)^{\frac {n-2}4}$,
conformal to $g$ and with $g$ a constant multiple
of $g_{\bar u}$.  For each $\e\in (0,\bar u]$ we have
$R(g_{\e})=n(n-1)$ and $g_{\e}$ is rotationally invariant with respect
to the $S^{n-1}$ factor and periodic in $t$.
Because of their similarity to the CMC surfaces of revolution
discovered by Delaunay \cite{D} these are called Delaunay solutions,
although it was Fowler \cite{F1}, \cite{F2} who first studied
the differential equation of which these are solutions.
Here $g_{\e}= u_{\e}^{\cl}g$, and we have normalized
so that at $t=0$, $u_{\e}^{\prime}=0$ and $u_{\e}^{\prime\prime}\ge  0$.
In general, these metrics also have a translation parameter which is
relevant to the analysis, as will be apparent below.
These solutions are discussed at length in \cite{MPU},
to which we refer for details on the discussion below,
cf. also \cite{S2}.

As $\e \rightarrow 0$, the
supremum of $u_\e$ is uniformly bounded, but the infimum tends to zero.
Geometrically, the metrics $g_{\e}$ develop a sequence of evenly spaced
`necks' which separate almost spherical regions. As $\e \rightarrow 0$,
these metrics converge to a `string of pearls' -- a sequence of round
spheres of radius one adjoined at their poles and arranged along a
fixed axis. We will denote $(\RR \times S^{n-1}, g_\e)$ by $D_\e$.

The linearization $\Le = \Lap_\e + n$ at any $g_{\e}$ is self-adjoint.
It has periodic coefficients, hence its spectrum is pure absolutely
continuous; there is no point spectrum, i.e. no eigenfunctions in $L^2$.
This last assertion may be seen rather concretely. Separating variables
according to the eigenfunction decomposition of $\Lap_\th$, we reduce
to analyzing each of the ordinary differential operators $L_j$ induced
on the eigenspaces. When $j > 0$, i.e. when the eigenfunction $\psi_j(\th)$
on $S^{n-1}$ is nonconstant, any solution of $L_j \phi = 0$ grows
exponentially in one direction or the other, as may be determined
by simple estimates \cite{MPU}. On the other hand, the two linearly
independent solutions of $L_0 \phi = 0$ are temperate, and so
we must examine them further to ensure that they are not in $L^2$.
Fortunately they can be determined explicitly by differentiating
$u_{\e}$ with respect to either the translation parameter $t$, or the
Delaunay parameter $\e$. Calling these $\phi_0^+$
and $\phi_0^-$, respectively, then $\phi_0^+$ is periodic, hence
bounded, while $\phi_0^-$ grows linearly in $t$.

It is crucial in what follows that these temperate Jacobi fields
are integrable, i.e. they arise as derivatives of one parameter
families of conformally related CPSC metrics.

Although $\Le$ does not have closed range on $L^2$, it does
have this property when considered as an operator on certain
weighted Sobolev or H\"older spaces. This was proved for the
weighted Sobolev spaces $H^s_\d$ in \cite{MPU}, and for
the weighted H\"older spaces in \cite{MPa2}. There are
advantages and disadvantages in using either of these spaces:
the H\"older spaces are better suited to the nonlinearity, but for
the various duality arguments we use the Sobolev spaces are more
convenient. Thus, for $\d, s \in \RR$, define
$$
H^s_\d(D_\e) = \{\phi = e^{\d\sqrt{1+t^2}} \tilde\phi:
\tilde \phi \in H^s(D_\e) \}, \tag 2.4
$$
where $H^s$ is the standard (global) Sobolev space on $D_\e$
with respect to $g_\e$, and $t$ is the `cylindrical length'
coordinate on $\RR \times S^{n-1}$. Note that the geometric length
along $D_\e$ depends on $\e$; for each $\e>0$ it is commensurate with
$t$, but not uniformly so as $\e \rightarrow 0$.
The spaces $H^s$ are algebras provided $s > n/2$, and so we
shall always make this restriction  when needed for the nonlinear
aspects of our problem.

Now
$$
\Le: H^{s+2}_{\d}(D_\e) \longrightarrow H^s_{\d}(D_\e) \tag 2.5
$$
is bounded for any $s$ and $\d$, but when $\d = 0$ it does not
have closed range. In fact, $\Le \phi = 0$ has a two dimensional
family of temperate solutions
(namely the span of $\phi_0^+$ and $\phi_0^-$)
and these may be used to construct an
orthonormal sequence $\phi^{(j)}$ in $H^0_{0}=L^2$ (or any
$H^s_0$) with $\|\Le \phi^{(j)}\| \rightarrow 0$; this is one of
the standard criteria for showing the range is not closed.

We prove in \cite{MPU} that there is a monotone sequence $\delta_j
\rightarrow \infty$ , depending on $\e$ and with $\delta_0 = 0$,
for which the map (2.5) is Fredholm provided $\d \notin \{\pm \delta_j\}$.
The values of the $\delta_j$ are exactly those for which the ordinary
differential operators $L_j$ have solutions $L_j \phi = 0$ growing exactly
like $e^{\pm(\d/P_{\e}) t}$, where $P_{\e}$ is the period of the metric
$g_\e$. (This is demonstrated in Proposition 4.8 of \cite{MPU}
under the assumption that $P_\e = 1$. A simple rescaling leads
to this version.)
Thus the same argument as above shows that $L_j$, and hence $\Le$, cannot
have closed range on $H^s_{\pm \d_j}$. The main content of this result
is that $\Le$ has closed range when $\d \ne \pm \d_j$.

There is no solution of $L_j \phi = 0$ which decays faster
than $e^{-|t|\d}$, for any fixed $\d > 0$, as $t \rightarrow \pm \infty$; for
$j \ge 1$ this follows from the maximum principle (cf.
\cite{MPU} for the case $j=1$ which is somewhat more subtle),
while for $j=0$ it follows because we know the solutions explicitly.
This implies that (2.5) is injective provided $\d < 0$.
By duality and elliptic regularity, (2.5) has dense range if $\d > 0$;
when $\d > 0$ and $\d \ne \d_j$, the range is also closed, hence (2.5)
is then surjective.

We have established that $\Le$ is surjective on $H^{s+2}_\d$ for
$\d > 0$, $\d \ne \d_1, \d_2, \dots$. Unfortunately, none of these
spaces are suitable for the nonlinear problem: if $\phi$ grows like
$e^{t\d}$ then $(1+\phi)^\cri$ grows even faster. It is possible
to obtain surjectivity on a smaller space. Let $\chi$ be a cutoff function
which equals one for $t \ge 1$ and zero for $t \le -1$, and define
the `deficiency subspace' $W$ by
$$
W = \text{span\,}\{\chi \phi_0^+, \chi \phi_0^-\}. \tag 2.6
$$

\proclaim{(2.7) Proposition} The map
$$
\Le: H^{s+2}_{-\d} (D_\e) \oplus W \longrightarrow H^s_{-\d}(D_\e)
$$
is surjective for any $\d < \d_1$.
\endproclaim

An analogous result is proved in \cite{MPU} for more general manifolds
of CPSC with $k$ asymptotically Delaunay ends, and with respect to the
weighted Sobolev spaces, but is quite simple to prove for the $D_\e$.
Suppose $f \in H^{s+2}_{-\d}(D_\e)$, and
let $f_j$ be its eigencomponents with respect to the Laplacian on $S^{n-1}$.
A solution $u_j$ of $L_j u_j = f_j$
may be constructed for each $j$ by `integrating in from $-\infty$' in
the standard ODE variation of parameters formula. These solutions
$u_j$ will decay like $e^{-|t|\d}$ as $t \rightarrow
-\infty$. Since $\d < \d_1$, $u_j$ must also decay like $e^{-|t|\d}$
as $t \rightarrow +\infty$ for $j \ge 1$, for if it did not,
then the difference between this solution and any other solution $v_j$
of $L_j v_j = f_j$ would be in the nullspace of $L_j$, hence not in
$H^{s+2}_\d$. This argument fails for $j = 0$, so $u_j$ can be written,
for $t \gg 0$, as a sum of two terms, one in $H^{s+2}_{-\d}$ and one in $W$.

We still need to check that the nonlinear operator $N_{g_\e}$
maps $H^{s+2}_{-\d} \oplus W$ to $H^{s}_{-\d}$. Clearly
$N_{g_\e}$ maps $H^{s+2}_{-\d}$ to $H^{s}_{-\d}$. To ensure that it
also carries $W$ to $H^{s}_{-\d}$ we need to modify the
definition of this map slightly. In fact, since elements
$(a\chi\phi_0^+, b\chi\phi_0^-) \in W$ correspond,
for $t \ge 1$, to the infinitesimal variations of one-parameter families
of Delaunay metrics, we can define a two parameter family of metrics
$\tilde{g}_{\e,a,b}$ on the cylinder such that for $t \le -1$,
$\tilde{g}_{\e,a,b}= g_{\e}$ and for $t \ge 1$,
$\tilde{g}_{\e,a,b}= g_{\e+d_{\e}(b)}(t-\tau_{\e}(a))$.
Here $d_{\e}:\RR\rightarrow (-\e,{\bar u}-\e)$  and
$\tau_{\e}:\RR\rightarrow (-P_{\e}/2,P_{\e}/2)$ are
monotone, smooth, surjective functions such that
$d_{\e}(0)=\tau_{\e}(0)=0$.
The map $(a,b)\mapsto \tilde{g}_{\e,a,b}$
induced by $\tau_{\e}$ and $d_{\e}$ can be
regarded as an exponential map to the space of Delaunay metrics on
the half cylinder from the tangent plane at the point $g_\e$.
By judicious choices of
the functions $d_{\e}$  and $\tau_{\e}$ and  the definition of
$\tilde{g}_{\e,a,b}$ in $-1<t<1$, we can insure that if
$\tilde{g}_{\e,a,b}=\tilde{u}_{\e,a,b}^{4/(n-2)}g_{\e}$, then
$$
\aligned
\chi\phi_0^+ = \chi\phi_0^{+}(\e)& = \frac {d}{da} \left.
\tilde{u}_{\e,a,b}\right|_{(a,b)=(0,0)} \\
\chi\phi_0^- = \chi\phi_0^{-}(\e)& = \frac {d}{db} \left.
\tilde{u}_{\e,a,b}\right|_{(a,b)=(0,0)}.
\endaligned
$$
We define a new operator by
$$
N_{g_\e}^{(a,b)}(\phi) = \Lap_{\tilde{g}_{\e,a,b}} \phi - \frac{n-2}{4(n-1)}
R(\tilde{g}_{\e,a,b})(1+\phi) + \frac{n(n-2)}{4}(1+\phi)^{\cri} \tag 2.8
$$
for $\phi \in H^{s+2}_{-\d}$ and $(a,b)\in\RR^2$.  Finally, setting
$$
N_{g_\e}(\phi,a\chi\phi_0^+, b\chi\phi_0^-)\equiv N_{g_\e}^{(a,b)}(\phi),
\tag 2.9
$$
we see that $N_{g_\e}:H^{s+2}_{-\d} \oplus W \rightarrow H^{s}_{-\d}$
is a well defined real analytic map.

\head Complete CPSC metrics on $X \setminus \Lambda$ \endhead

In the last subsection we considered singular Yamabe metrics
on the sphere with two points removed. Rather different
solutions were constructed in \cite{MPa1} and \cite{MS}.
These are complete CPSC metrics on $M = X \setminus \La$ where $\La$ is a
finite disjoint union of submanifolds $\La_i$ without boundary,
$(X,g_0)$ is compact of nonnegative scalar curvature, and
$\dim\La_i = k_i$ with $1 \le k_i \le (n-2)/2$. The upper bound
on the dimensions $k_i$ is, by a theorem
of Schoen and Yau  \cite{SY}, a necessary condition.
Note that we temporarily abandon the convention that $g_0$ is
complete on $M= X\setminus\La$ here.

The completeness of $g = u^{\cl}g_0$ on $M\setminus \La$ necessitates that
$u$ tends to infinity rather strongly on approach to $\La$. A detailed
study of this singular behavior is given in \cite{M1}. Let $r$ denote
a smooth function on $M$ which is everywhere positive, and which agrees
with the polar distance function (with respect to $g_0$) on a tubular
neighbourhood of $\La$. The solutions $u$ constructed in \cite{MPa1}
and \cite{MS} are asymptotic, to leading order, to $Ar^{(2-n)/2}$
as $r \rightarrow 0$, with $A$ a constant
depending only on dimension. It is shown in \cite{M1} that in a
neighbourhood of each component $\La_i$ these solutions have more
refined asymptotics $u \sim Ar^{(2-n)/2}(1 +  O(r^{k_i/2}))$.
For convenience in the rest of this section,
assume that $\La$ has only one component of dimension $k$.

The linearized scalar curvature operator relative to one of these
CPSC metrics $g$ has the form (2.3). If $k^2 \le 4(n-2)(n-2k-2)$
(cf. \cite{MS}), the continuous spectrum of $L$ contains $0$, hence
again $L$ does not have closed range on $L^2$.  As before, it is
appropriate to let $L$ act on weighted Sobolev or H\"older spaces.
Although both \cite{MPa1} and \cite{MS} use H\"older spaces, we
shall use Sobolev spaces as above. Once again, $H^s_\d(M,g)$
is defined to be the space of functions $\phi = r^{-\delta + \frac{k}{2}}
\tilde\phi$, where $\tilde \phi$ is in the uniform global Sobolev space
$H^{s}$ on $M$ with respect to the complete metric $r^{-2}g_0$ (or,
equivalently, with respect to $g$). Note the change of sign and
shift of weight parameter relative to the previous definition.
Then it follows from the theory of \cite{M2} that
$$
L: H^{s+2}_{\d}(M) \longrightarrow H^{s}_{\d}(M)
\tag 2.10
$$
has closed range for all $\d \notin \{\pm \delta_j\}$, where
as before $0 = \d_0 < \d_1 < \dots \rightarrow \infty$.
Unlike the situation for the Delaunay metrics, (2.10) will not
be Fredholm, even when it has closed range. Indeed, for $\d < 0$ it
has infinite dimensional kernel, but at most finite dimensional cokernel,
while for $\d > 0$ it has infinite dimensional cokernel and at most
finite dimensional kernel.

The CPSC metrics constructed in \cite{MS} and \cite{MPa1} are
nondegenerate in the sense that (2.10) is surjective
if $0 = \d_0 < \d < \d_1$ (and hence for all $\d > 0$, $\d \ne \d_j$).
As shown in \cite{MS}, this implies that every sufficiently
small element of the nullspace of $L$ in $H^{s+2}_{\d}$ for
$0 < \d < \d_1$ is integrable, i.e. is the tangent vector of a one-parameter
family of solutions $u_t$.  Notice that because of the shift
in the weight parameter here, the space on which $L$ is surjective contains
only decaying functions, so unlike before we do not need to separate off the
nullspace (or deficiency subspace) as in (2.7) to obtain a space on which
the nonlinear operator $N_g$ acts.

\head CPSC metrics on manifolds with boundary\endhead

Our construction also applies to CPSC manifolds with boundary, either
compact or noncompact.  The issue here is the mapping properties
of $L$ on Sobolev or H\"older spaces with Dirichlet boundary values.
A geometrically natural boundary condition for the nonlinear problem
is to require the boundary in the induced metric to have constant mean
curvature. This has been studied extensively by Escobar \cite{E} and others.
Two key examples are the spherical cap $S^{n}_{r}$ of radius $r$
and the Delaunay metrics on the half cylinder, $D_{\e}^{\alpha}$, which
are simply the restrictions of the $D_\e$ to $t \ge \alpha$.
The mean curvature of the boundary is a constant depending on
$r$ and $\alpha$; when $r = \pi/2$, $\alpha = 0$ the boundaries
are not only minimal but totally geodesic.

The half-Delaunay metrics are all nondegenerate, as defined below.
This follows from simple modifications of the previous discussion
of the full Delaunay metrics. On the other hand, the spherical
cap $S^n_r$ is nondegenerate only when $r \ne \pi/2$. When $r = \pi/2$,
$L$ has a one-dimensional nullspace consisting of the linear functions
vanishing on the boundary, hence by duality is not surjective.

\head Nondegeneracy \endhead

Having described in some detail the main examples of CPSC manifolds
for which the Jacobi operators $L$ are in some sense surjective,
we now abstract these properties and formulate a general notion of
nondegeneracy sufficiently flexible for the gluing construction.

Suppose $(M,g)$ is a noncompact, complete Riemannian manifold
of CPSC. The standard Sobolev spaces $H^s$ are defined relative
to the Riemannian measure and connection. We shall assume that
there exists a weight function $0 < \a \in \Cal C^\infty(M)$ the
powers of which define a scale of weighted $L^2$ and Sobolev spaces.
Thus we define
$$
H^s_\d(M) = \{v = \a^\d \tilde v: \tilde v \in H^s(M)\}. \tag 2.11
$$
The dual of $H^s_\d$ is naturally identified with $H^{-s}_{-\d}$.

The main nondegeneracy hypothesis is that there exists a weight parameter
$\d>0$ such that for all $s \in \RR$ there exists a constant $C = C_s > 0$
for which
$$
\| \phi \|_{s+2,-\d} \le C \| L \phi \|_{s,-\d} \tag 2.12
$$
for every $\phi \in \cin_0(M)$. This implies that
$$
L: H^{s+2}_{-\d} \longrightarrow H^s_{-\d} \tag 2.13
$$
is injective and has closed range. It also gives some analytic control on the
behavior of $L$ on the ends of $M$. By duality we see that
$$
L:H^{s+2}_{\d} \rightarrow H^{s}_{\d} \tag 2.14
$$
is surjective. (It is precisely this last assertion which, though
still true, would be a bit more difficult to obtain if we were using
H\"older spaces.)

In some cases, such as for the problem on $M\setminus\La$ where
all components of $\La$ are of positive dimension, this is the
only hypothesis needed because, for some neighbourhood of zero
$\Cal U \subset H^{s+2}_\d$, and for $\d > 0$ sufficiently small
$$
N_g: H^{s+2}_\d \longrightarrow H^s_\d  \tag 2.15
$$
is well defined and has surjective linearization.
In other cases, such as for the Delaunay metrics, $N_g$ does
not map elements of $H^{s+2}_\d$ to $H^s_\d$ and so we need to
find another space on which the linearization is surjective
and on which $N_g$ is well-behaved. Thus we assume the existence
of a `deficiency space' $W \subset H^{s+2}_\d$, composed
of elements of the form $\chi \phi$, for some fixed cutoff function $\chi$,
where $\phi \in H^{s+2}_\d$ and $L\phi = 0$ outside some compact set,
such that
$$
L: H^{s+2}_{-\d} \oplus W \longrightarrow H^{s}_{-\d} \tag 2.16
$$
is surjective. There is no loss of generality in assuming that (2.16)
is an isomorphism, because we can always restrict to the orthogonal
complement of the intersection of the nullspace of $L$ on $H^{s+2}_\d$
with $H^{s+2}_{-\d} \oplus W$, no element of
which is contained in $H^{s+2}_{-\d}$ by hypothesis.
We also require that the elements of $W$ are `asymptotically
integrable,' which we take to mean that elements of $W$ are derivatives
of one parameter families of exact solutions of $N_g$ outside a compact set.
The validity of this condition was discussed in detail for the
Delaunay metrics. Rather than formulate the asymptotic integrability
more specifically, we refer to these examples and single out as the second
nondegeneracy hypothesis the only consequence that we require, namely that
$$
N_g: \Cal U \longrightarrow H^s_{-\d} \tag 2.17
$$
is surjective onto a neighbourhood of $0$ with surjective
linearization (2.16), where $\Cal U$
is a neighbourhood of the origin in $H^{s+2}_{-\d} \oplus W$.

It turns out that this second condition is rather less general than
it might appear. In fact, it is not hard to show that (2.16) implies
that $W$ must be finite dimensional. For if it were not, then one could
construct an orthonormal sequence $\{\chi \phi_j\}\subset W$, of fixed
$H^{s+2}_\d$ norm one, which decays to zero uniformly on any compact set.
This would contradict the closedness of the range of (2.16).

We say that a metric $g$ is nondegenerate if the linearization at $1$
of the nonlinear operator $N_g$ for $g$ is nondegenerate, in the sense that
(2.12) and either (2.15) or (2.16-2.17) hold.
Note that in the two main examples indicated above
nondegeneracy holds provided that $L$ has no kernel in $L^2$.
The analytic control at infinity which improves this to (2.12) and (2.16)
is provided by the strong asymptotics which these solutions exhibit,
cf. [MPU], [M1] and [M2].

\specialhead III. The gluing construction \endspecialhead

Our aim in this section is to state and prove the gluing
theorems for manifolds with nondegenerate CPSC metrics.
Thus let $(M_1,g_1)$ and $(M_2,g_2)$ be two complete CPSC manifolds
(possibly with boundary).  In the next subsection we will construct
a one-parameter family of approximate solution metrics $g_\e$
on the connected sum $M_1 \#_\e M_2$, where the parameter $\e$ corresponds
to the size of the connecting neck and is assumed to be small.
The approximate solution metric $g_{\e}$ has CPSC except on a neighbourhood
of this neck.

\proclaim{(3.1) Theorem} Suppose $(M_1,g_1)$ and $(M_2,g_2)$ are
two nondegenerate CPSC manifolds.
Then for some $\e_0 > 0$ and all $0 < \e < \e_0$,
there exists a function $u \in H^{s+2}_{-\d}(M_1 \#_\e M_2, g_\e)$
such that ${\bar g}_{\e}=(1 + u)^{\cl}g_\e$ is nondegenerate with CPSC.
\endproclaim

In the next section we will make more refined statements about the
global geometry of these new metrics and the implications of this
construction for the moduli spaces.

In the rest of this section we shall prove Theorem (3.1). The proof has
several steps. We first construct the approximate solution metrics $g_{\e}$.
In the next two steps, which are the heart of the proof, we show that
the $g_\e$ are nondegenerate and that (right) inverses for the linearized
scalar curvature operators are uniformly bounded as $\e \rightarrow 0$.
The rather simple indirect method  used here is the main novel
ingredient in this paper. Finally, using this
nondegeneracy, we perturb $g_\e$ to an exact solution using
a standard iteration argument.

\head Approximate solutions \endhead

Let $(M_i,g_i)$, $i = 1,2$, be two nondegenerate complete
CPSC manifolds. Fix points $p_i \in M_i$ and small metric balls
$B_{2\a_i}(p_i)$. Let $(r_i,\th_i)$ be Riemannian polar coordinates
about $p_i$. Then for each $\e \in (0,1)$ identify the
annulus $B_{\a_1}(p_1) \setminus B_{\e\a_1}(p_1)$ with
$B_{\a_2}(p_2) \setminus B_{\e\a_2}(p_2)$ by the relation $(r_1,\theta_1)
\sim (r_2, \theta_2)$ if $\theta_1 = \theta_2$ and $r_1 r_2 = \e \a_1 \a_2$.
This is the connected sum $M_{\e} \equiv M_1 \#_{\e} M_2$; the points
$p_i$ and radii $\a_i$ are suppressed in this notation, although both
metrically and conformally this data is important.

We first consider the case where the metrics $g_i$ are conformally
flat in the balls $B_{2\a_i}(p_i)$.
In this case the analysis is the most transparent.
It is useful to rephrase the problem relative
to new, conformally equivalent, background metrics $g_{i,c}$ on
$M_i \setminus p_i$. Geometrically we deform the conformally flat
metrics $g_i$ in $B_{\a_i}(p_i)$ to half-infinite cylinders.
The connected sum $M_\e$ is then given by identifying these
cylinders at a certain distance.
The metric degeneration of $M_\e$ as $\e\rightarrow 0$ now
corresponds to the lengthening  of this cylindrical tube.
More specifically, by (temporary) hypothesis, $g_i = v_i^{\cl}\delta$
in $B_{2\a_i}(p_i)$, where $\delta$ is the standard Euclidean
metric. We set $g_{i,c} = u_i^{-\cl} g_i$, where
$$
u_i = \rho_i + (1-\rho_i)
((n-2)/n)^{\frac{2-n}4}r^{\frac{n-2}2}v_{i},
$$
for some smooth cutoff function $\rho_i \ge 0$ with $\rho_i = 1$ on
$M_i\setminus B_{2\a_i}(p_i)$ and $\rho_i = 0$ in $B_{\a_i}(p_i)$. Then in
$B_{\a_i}(p_i) \setminus p_i$, $g_{i,c}$ is isometric to the
cylindrical metric $ \frac{n-2}{n}(dt_i^2 + d\th_i^2)$. The normalizing
constant is chosen so that this cylinder has scalar curvature $R = n(n-1)$.

In this smaller ball replace the variable $r_i$ by $t_i = -\log r_i$,
and let $T = T(\e) = -\log \e$. Then the identification between the two
annular regions is now given by $(t_1,\theta_1) \sim (t_2,\theta_2)$ if
$\th_1 = \th_2$ and $t_1 + t_2 = A_1 + A_2 + T$,
where $A_i = -\log \a_i$ and $A_i \le t_i \le A_i + T$.
We now alternately denote the connected sum $M_\e$ by $M_T$.
Because we have assumed that the $g_i$ are conformally flat on
$B_{2a_i}(p_i)$,
this identification map is an isometry with respect to the  metrics
$g_{i,c}$, hence there is a naturally induced metric $g_{c,T}$ on $M_T$.
We let $C_T$ denote the cylindrical region where $A_i \le t_i \le A_i + T$.

The approximate solution metric $g_T$ on $M_T$ is defined in
terms of the conformal factor
$$
u_T = \chi_1 u_1 + \chi_2 u_2,
$$
using nonnegative cutoff functions $\{\chi_1,\chi_2\}$ on $C_T$, where
$\chi_i \equiv 1$ for $t_i \le A_i + T/2 - 1$
and $\chi_i = 0$ for $t_i \ge A_i + T/2 + 1$ (here we regard $\chi_i$
as a function on $M_i$). $u_T$ extends naturally to all of $M_T$, and
we define the approximate solution metric by
$$
g_T = u_T^{\cl} g_{c,T}. \tag 3.2
$$
Note that $g_T = g_i$ on $M_i\setminus B_{c(\e)\a_i}(p_i)$, where
$c(\e)=c{\sqrt\e}$.

If the metrics $g_i$ are not conformally flat in the balls $B_{2\a_i}(p_i)$,
the construction of $g_T$ is almost identical, but is no
longer conformally natural, i.e. the conformal class of $g_T$
depends on choices of cutoff functions as well as $T$.
In $B_{2\a_i}(p_i)$ we can choose a normal coordinate system
in terms of which $g_i = \delta + h_i$, where $h_i$ is
small in some fixed norm. In $B_{2\a_i}(p_i)\setminus B_{\a_i}(p_i)$
we deform $h_i$ to zero and simultaneously deform $\delta$
to $((n-2)/n)r^{-2}\delta$. These metrics may now be joined
as before.

We define the unweighted Sobolev spaces $H^s$ with respect to
the metrics $g_{c,T}$; of course, the norms on these spaces
depend on $T$, although this effect may be localized to $C_T$.
We may assume that the weight functions $\a_i$ on $M_i$ are
identically one in a large neighbourhood of the points $p_i$;
these extend naturally over $C_T$ and we may define a new
weight function $\a$ on $M_T$. Using it, we define the
weighted Sobolev spaces $H^s_\d$ on the connected sum.

The metric $g_T$ does not have CPSC, and we can easily estimate
the error term
$$
f_T = \Lap_{g_{c,T}} u_T - \frac{n-2}{4(n-2)}R(g_{c,T})u_T + \frac{n(n-2)}{4}
u_T^{\cri}.  \tag 3.3
$$
In fact, (in the conformally flat case) $g_{T}$ does have scalar
curvature $n(n-1)$ except in the middle of $C_T$, where $t_i \in
[A_i + T/2-1, A_i + T/2 + 1]$. Since $u_i = O(r_i^{\frac{n-2}{2}})$
and $r_i = e^{-t_i}$, it is clear that
$$
\|f_T\|_{s} \le C e^{-T/2} \tag 3.4
$$
for any $s$. We do not need to use a weighted norm here because
$f_T$ is supported on $C_T$ where $\a = 1$.
In the non conformally flat case, there is an additional error
term incurred by cutting off $h_i$. By choosing the normal
coordinates correctly, we can also bound this extra error term, which is
now supported near the boundary of $C_T$, by $Ce^{-T/2}$.

\head Nondegeneracy of the approximate solution \endhead

In order to be able to perturb $g_T$ to a CPSC metric, we need
to establish nondegeneracy of the linearization of the scalar curvature
operator for  $g_{c,T}$ at $u_T$. Although the definitions
(2.12) and (2.15-17) of nondegeneracy were given only for CPSC metrics,
these hypotheses make perfect sense here.

\proclaim{(3.5) Proposition} There exists a $T_0 > 0$ such that for all
$T \ge T_0$, the metric $g_T$ of (3.2) is nondegenerate.
\endproclaim

Before embarking on the proof, we make some preliminary observations
about the linearizations of the scalar curvature operator on
$M_1$, $M_2$ and $M_T$.  For any metric $g$, the conformal Laplacian
$$
\Cal L_g = \Lap_g - \co R(g_0)
$$
is the linear part of $N_g$ in (2.1). It is conformally
equivariant in the sense that if $g' = u^\cl g$,
then for any $\phi \in \Cal C^{\infty}$,
$$
\Cal L_g (u\phi) = u^\cri \Cal L_{g'} \phi. \tag 3.5
$$
A special case of this equality is when $\phi = 1$, in which
case (3.5) reduces to (2.1).

Next, suppose $g' = u^{\cl}g$, and let $L_g$ be the linearization
of $N_g$ at $u$.  The relationship between $L_g$ and $L_{g'}$
is not as simple as (3.5) in general, but it is when both $g$ and $g'$ have
the same (constant) scalar curvature. Indeed, if $R(g) = R(g') = n(n-1)$,
then  $\Cal L = \Lap - \frac{n(n-2)}{4}$ for $g$ or $g'$. Since
$L_g = \Cal L_g + \frac{n(n+2)}{4} u^{\cl}$, we have
$$
\gathered
L_g( u\phi)  = \Cal L_g(u\phi) + \frac{n(n+2)}{4}u^\cl (u\phi) \\
= u^{\cri}\left( \Cal L_{g'} \phi + \frac{n(n+2)}{4} \phi \right)
= u^{\cri} L_{g'}\phi .
\endgathered
\tag 3.6
$$
In particular, away from the transition regions $B_{2\a_i}(p_i)\setminus
B_{\a_i}(p_i)$, (3.6) applies to either of the two pairs of
metrics $g_i$, $g_{i,c}$, with $u = u_i$. The linearizations
corresponding to these two metrics will be denoted
$L_i$ and $L_{i,c}$.

\proclaim{(3.7) Lemma} Suppose $L_{i,c} \phi = 0$ for
some function $\phi$ on $M_i\setminus \{p_i\}$, and suppose
that $\phi$ is bounded on the deleted neighbourhood
$B_{\a_i}(p_i)\setminus \{p_i\}$. Then $u_i^{-1} \phi$
extends smoothly across $p_i$ on $M_i$, and in particular,
$|\phi| \le Ce^{-(n-2)t_i/2} = Cr_i^{\frac{n-2}{2}}$ on this neighbourhood.
\endproclaim

\demo{Proof:} By (3.6)
$$
u_i^{\cri}L_i (u_i^{-1} \phi) = L_{i,c} \phi = 0,
$$
and so, letting $\psi = u_i^{-1} \phi$, we see that $L_i \psi = 0$
on $B_{\a_i}(p_i)\setminus \{p_i\}$. Since $\phi$ is bounded,
$|\psi | \le Cr^{(2-n)/2}$. Thus $\psi$ extends to a weak solution of
$L_i\psi = 0$ on all of $M_i$, and by a standard removable
singularities theorem extends smoothly across $p_i$. \qed
\enddemo

One other result is needed. Let $L_T$ denote the linearization
of $N_{g_{c,T}}$ at $u_T$.
\proclaim{(3.8) Lemma} Suppose that $\phi$ solves $L_T \phi = 0$
on the cylindrical region $C_T$, and furthermore suppose that
$\|\phi\|_{L^2(A)} \le 1$, where $A$ is the union of two annular
neighbourhoods, one about each of the boundary components
of $C_T$. Then $|\phi| \le C$ on all of $C_T$, where $C$ is
independent of $T$.
\endproclaim
\demo{Proof} Since $g_{c,T}$ is a product metric in $C_T$,
$$
L_T = \frac{n}{n-2}\left( \frac{\del^2\,}{\del t^2} + \Lap_{\th}\right)
- \frac{n(n-2)}{4} + \frac{n(n+2)}{4}u_T^{\frac{n+2}{n-2}}
$$
there. Provided we adjust the annular region appropriately, we can
ensure that the term of order zero in this operator is strictly
negative. The result then follows from the maximum principle, and
in fact we can take $C = 1$. \qed
\enddemo
\flushpar It is not strictly necessary that the term of order zero
is negative on the whole cylinder; it is only necessary that it
is nonnegative on a compact set not growing in size as $T$ gets
large. We leave details to the reader.

\demo{Proof of Proposition (3.5)} We first show that the linearization
$L_T$ of $N_{g_{c,T}}$ at $u_T$  is injective on $H^{s+2}_{-\d}$ for any
$s$, and we argue by contradiction.  Suppose that there
exists a sequence $T_j \rightarrow \infty$ and a function $\phi_j \in
H^{s+2}_{-\d}(M_{T_{j}})$ such that $L_{T_j}\phi_j = 0$. The weight $-\d$ of
course refers to growth behavior on any other ends of $M_1$ and $M_2$.
Choose compact neighbourhoods $K_i \subset M_i \setminus \{p_i\}$
containing $\del B_{\a_i}(p_i)$ and normalize $\phi_j$ so that
$$
\max_{i=1,2}\left\{ \|\phi_j\|_{H^{s+2}(K_1)},
\|\phi_j\|_{H^{s+2}(K_2)} \right\} = 1.
$$
In particular, on one or the other of these sets $\phi_j$ has
Sobolev norm uniformly bounded below; by passing to a subsequence we
can assume that this takes place on $K_1$. Since we can assume that $s >
n/2$, by elliptic regularity, we also get uniform supremum bounds for
$\phi_j$ on $K_1$ and $K_2$.

We can now take the limit as $j \rightarrow \infty$
to obtain a limit $\phi$, which is a function on the disjoint
union $M_1\setminus\{p_1\} \sqcup M_2 \setminus \{p_2\}$.
By the uniform lower bound, $\phi$ is nontrivial on $K_1$ and
solves $L_{1,c} \phi = 0$ there. Using Lemma (3.8) we see that $\phi$
is bounded along the cylindrical end of $M_1 \setminus \{p_1\}$,
hence, by Lemma (3.7), the function $\psi = u_1^{-1} \phi$ extends
smoothly to all of $M_1$ and satisfies $L_1 \psi = 0$.
It is also the case that $\psi \in H^{s+2}_{-\d}(M_1)$. To see this let
$\chi$ be a cutoff function equaling $1$ outside of $B_{2\a_1}(p_1)$,
vanishing near $p_1$, and with the support of $\nabla\chi\subset
K_1$.
Then
$$
L_T (\chi \phi_j) = L_{1,c}(\chi \phi_j) = (\Lap\chi)\phi_j +
2 \nabla \chi \cdot \nabla \phi_j.
$$
The right hand side is compactly supported and uniformly bounded
in $H^s_{-\d}(M_1)$. By (2.12) we obtain a uniform bound for
$\chi \phi_j$ in $H^{s+2}_{-\d}$, and since $u_1 = 1$ outside
$B_{2\a_1}(p_1)$ it is easy to see that $\psi \in H^{s=2}_{-\d}(M_1)$,
as claimed. This is a contradiction, since $\psi$ is nontrivial
and $L_1$ satisfies (2.12). This proves the injectivity of
$L_T$ for $T$ sufficiently large. Now we may patch together the
estimates (2.12) from $M_1$ and $M_2$ to obtain, for each $T$,
$$
\| \phi \|_{s+2,-\d} \le C( \|L_T \phi \|_{s,-\d} + \| \phi \|_{0,K})
$$
where the final term is the $L^2$ norm of $\phi$ on a fixed compact
set $K$. Finally, it is standard that the injectivity of $L_T$ shows
that this final term must be bounded by a fixed constant multiple
(possibly depending on $T$) of $\|L_T \phi\|_{s+2,-\d}$, so we have shown
that $(M_T,g_T)$ satisfies (2.12) for $T$ sufficiently large.

If we are in a case where (2.15) holds, then we have proved that
$g_T$ is nondegenerate already. Thus suppose that either $M_1$ or $M_2$,
or both, have deficiency spaces $W_1$ and $W_2$ satisfying (2.16).
Let $W = W_1 \oplus W_2$. We must show that $L_T$ satisfies
(2.16) for $T$ sufficiently large.  For clarity, we denote the
restricted (or extended, depending on your viewpoint) map in (2.16)
by $\tilde L_T$. We shall verify that $\tilde L_T$ is surjective
by computing its adjoint $\tilde L_T^*$ and showing that
it has no nullspace. Note that $\tilde L_T$ obviously has
closed range because it is a finite dimensional extension of
an operator with closed range.

We formulate this somewhat more generally in the

\proclaim{(3.9) Lemma} Assuming the general setup of \S 2 above,
the map $L$ in (2.16) is surjective
if and only if for every $\phi \in \Cal B \equiv \ker (L) \cap H^{s}_{\d}$,
the linear functional on $W$ defined by
$$
\int (Lw) \phi \tag 3.10
$$
for $w \in W$, is not identically zero.
\endproclaim

\demo{Proof} The surjectivity of the operator in (2.16), which we
denote by $\tilde L$ temporarily, is equivalent to the
injectivity its adjoint $\tilde L^*$. We first compute this adjoint.
If $(v,w) \in H^{s+2}_{-\d} \oplus W$ and $f \in H^s_{-\d}$,
then $\tilde L^*$ is defined by
$$
\int {\tilde L}(v+w)f \a^{2\d} = \int v (\a^{-2\d}L \a^{2\d} f) \a^{2\d}
+ \int L(w) (\a^{2\d}f).
$$
Note that we cannot integrate the second term by parts because
$w$ does not decay, even though $Lw$ is compactly supported.
Now multiplication by $\a^{2\d}$ defines the natural isomorphism
between $H^s_{-\d}$ and $H^s_\d$; since $L$ is self-adjoint
when $\d = 0$ we see that $\a^{-2\d}L \a^{2\d}$ is canonically
identified with $L$ on $H^s_\d$.

Now if $f$ is in the nullspace of $\tilde L^*$, then setting $w = 0$ and
letting $v$ range over all of $H^{s+2}_{-\d}$ we see that
$\a^{-2\d}L \a^{2\d}f = 0$. Equivalently, $\phi = \a^{2\d}f$
is in the nullspace of $L$ in $H^{s}_{\d}$. Now letting $v=0$ we see
that $\phi$ is orthogonal (in unweighted $L^2$) to $Lw$ for every $w \in W$.

On the other hand, if $\phi \in H^{s}_{\d}$ with $L\phi = 0$
and $\int (Lw)\phi = 0$ for every $w \in W$ then $f = \a^{-2\d}\phi$
is in the nullspace of $\tilde L^*$.  \qed
\enddemo

We return to the proof of the proposition. We need to show
that $\tilde L_T^*$ is injective for $T$ sufficiently large.
As usual, assume not, so that there exists a sequence $T_j$
tending to infinity and corresponding elements $\phi_j$
in the nullspace $\Cal B$ such that $\int \phi_j L_Tw = 0$ for
all $w \in W$. (Note that the space $W$ is independent of $T$ because
its elements are supported away from $C_T$.)
Normalize the sequence as before, so that its norm on one of the
two compact sets $K_i \subset M_i$ is one. By the finite dimensionality
of $W$ the restriction (of a subsequence) of the $\phi_j$ converges
on $M_1 \setminus \{p_1\}$, say, to a nontrivial element $\phi$. As before,
$\psi = u_1^{-1}\phi$ extends smoothly to all of $M_1$ and solves
$L_1 \psi = 0$. This is a contradiction since $\int \psi L_1 w_1 = 0$
for all $w_1 \in W_1$ implies by Lemma (3.9) and
the surjectivity of $L_1$ in (2.16) on $M_1$, that $\psi = 0$.  \qed
\enddemo

\head Uniform surjectivity of $L_T$ \endhead

By Proposition (3.5) and duality, $L_T$ is surjective on $H^s_{\d}(M_T)$
for all $T \ge T_0$. This means that there exists some right inverse
$$
G_T: H^s_\d(M_T) \longrightarrow H^{s+2}_\d(M_T). \tag 3.11
$$
Because $L_T$ is not injective on $H^{s+2}_\d$, there are many
choices for the map (3.11). The canonical choice is the one
which has range agreeing with the range of the adjoint map
$L_T^*$. Henceforth we assume that the map (3.11) satisfies this condition.

To understand this choice better, note that by what we have proved,
$L_T L_T^*$ is an isomorphism, hence has a unique inverse $\Cal G_T$.
Since $L_T L_T^* \Cal G_T = I$, we see that $G_T$ must be $L_T^* \Cal G_T$
because both are right inverses of $L_T$ with range contained in the range
of $L_T^*$. We note also that by the same formalism as discussed
in the last subsection using the weight function $\a$, we may
identify this adjoint $L_T^*$ with $\a^{2\d} L_T \a^{-2\d}$.

Now we may restrict $G_T$ to $H^{s}_{-\d}$. Unfortunately,
its range may not coincide with $H^{s+2}_{-\d} \oplus W$.
Thus for a given $f \in H^s_{-\d}$ we have found two
possibly distinct solutions of $Lu =f$, namely the
solution $v + w \in H^{s+2}_{-\d} \oplus W$ and $G_T f$.
The difference between these solutions is an element $\phi \in \Cal B$,
the nullspace of $L_T$ in $H^{s+2}_\d$.  Our ultimate
goal is to show that $\|v+w \|$ is bounded by a multiple of
$\|f\|$, uniformly in $T$. We do this in two steps, first
showing that the norm of $G_T$, and then the correction term $\phi$,
are uniformly bounded.

\proclaim{(3.12) Proposition} The norm of the map $G_T$ in
(3.11) is uniformly bounded as $T \rightarrow \infty$.
\endproclaim
\demo{Proof} The proof, once again, is indirect. Thus we assume that
the result is false, so that for some sequence $T_j \rightarrow \infty$
there are functions $f_j \in H^s_\d$ with $\|f_j\|_{s,\d} \rightarrow 0$,
such that $\|G_{T_j}f_j\|_{s+2,\d} = 1$. Since each $G_{T_j} f_j = \psi_j$
is in the range of $L_{T_j}^*$, there exist functions $v_j \in
H^{s+4}_{\d}$ with $L_{T_j}^*v_j = \psi_j$. Because $\|\psi_j \| = 1$
and because of the boundedness of the $L_{T_j}^*$ on Sobolev spaces,
we know that
$$
\| v_j \|_{s+4,\d} \ge C.
$$

Our goal is to show that some subsequence of the $v_j$ converges to
function $v \in H^{s+4}_\d$ which is nontrivial on at least one of
$M_1$ or $M_2$. Suppose for definiteness that $v \ne 0$ on $M_1$.
Because of the boundedness of $L_{T_j}^*$ on these Sobolev spaces, we
can also assume that $\psi_j$ converges to $\psi \in H^{s+2}_{\d}$
where $L_{1,c}^*v = \psi$; also $\|f_j\|_{s,\d} \rightarrow 0$, so
that $L_{1,c}\psi = 0$ and hence $L_{1,c} L_{1,c}^* v = 0$.
Using Lemma (3.7) in the same way as before we see that
$u_1^{-1}\psi \equiv \phi$ is smooth on $M_1$, so that
$|\psi| \le Ce^{(2-n)t/2}$ on the cylindrical end of $M_1$.
This allows us to integrate by parts to conclude that
$$
0 = \langle v, L_{1,c} L_{1,c}^* v \rangle = \langle \psi, \psi \rangle
$$
and so $\psi = 0$. But now $L_{1,c}^*v = 0$, or equivalently,
$L_1^*(u_1^{-1} v) = 0$. But $u_1^{-1} v \in H^{s+4}_\d(M_1)$
and we have already established that $L_1^*$ is injective on this
space (or rather, we established that $L_1$ is injective on
$H^{s+4}_{-\d}$, which is the same). This is a contradiction, hence
the $v_j$ cannot converge as claimed, and the maps $G_{T}$ must
be uniformly bounded.

To finish the proof we must show that the $v_j$ converge in $H^{s+4}_\d$.
First we transform the problem so that we may work in $H^s_{-\d}$
and so avail ourselves of the estimate (2.12). Let $\tilde v_j
= \a^{-2\d}v_j \in H^{s+4}_{-\d}$, and define $\tilde \psi_j$ similarly.
Then $L_{T_j} \tilde v_j = \tilde \psi_j$.  Let $\chi$ be a smooth,
nonnegative cutoff function on either $M_1$ or $M_2$ which equals one
outside $B_{2\a_i}(p_i)$ and zero inside $B_{\a_i}(p_i)$. By computing
$L_{T_j}(\chi \tilde v_j)$ we see that
$$
\|\chi \tilde v_j \|_{s+4,-\d} \le C \left( \|\chi \tilde \psi_j \|_{
s+2,-\d} + \|\tilde v_j \|_{s+3,K}\right),
$$
where $K$ is some compact set containing the support of $\nabla \chi$
and $\| \cdot \|_{s,K}$ is the Sobolev $H^s$ norm on $K$.

We now show that $\|\tilde v_j \|_{s+3,K}$ is bounded away from
zero and infinity as $j \rightarrow \infty$. The upper bound
will imply, by (2.12), that $\|\tilde v_j\|_{s+4,-\d} \le C$.
The lower bound will imply that $\tilde v_j$ does not converge to zero on
$K$. This will show that $\tilde v_j$, hence $v_j$ also, must converge
to a nonzero function, which we know from above cannot happen.

There are two cases. Suppose first that $\| \tilde v_j \|_{s+3,K} \rightarrow
\infty$. Rescale $\tilde v_j$ and $\tilde \psi_j$ by the factor
$\|\tilde v_j\|_{s+3,K}^{-1}$ to obtain functions $\bar v_j$,
$\bar \psi_j$ with $L_{T_j} \bar v_j = \bar \psi_j$,
and $\|\bar v_j \|_{s+3,K} = 1$, $\|\bar \psi_j\|_{s+2,-\d}
\rightarrow 0$. Since $\bar v_j$ has fixed norm on a fixed
compact set, $\|\chi \bar v_j\|_{s+4,-\d} \le C$ by (2.12), so
we may pass to a limit and obtain a function $\bar v \in H^{s+4}_{-\d}$
with $L_{j,c} \bar v = 0$ for $j = 1,2$. Furthermore,
the restriction of $\bar v$ to at least one of the $M_i$,
say $M_1$, is nonzero. The boundedness of $\bar v_j$ on $K$, which contains
a neighbourhood of the boundary of $C_{T_j}$, implies by Lemma~(3.8)
that $\bar v_j$, hence $\bar v$ too, are uniformly bounded on $C_{T_j}$.
As before, this shows that $\bar w = u_1^{-1} \bar v \in H^{s+4}_{-\d}(M_1)$.
But $L_1 \bar w = 0$, which is a contradiction since $\bar w \ne 0$.

The other case is that $\|\tilde v_j \|_{s+3,K} \rightarrow 0$. Use the
same cutoff function $\chi$ as above (say on $M_1$), to compute that
$$
L_{T_j} L_{T_j}^* (\chi v_j) =
\chi f_j + [L_{T_j}L_{T_j}^*, \chi] v_j \equiv h_j.
$$
By hypothesis, both terms on the right tend to zero in $H^s_{\d}$.
But $L_{T_j} = L_1$ on the support of $\chi$, and if the $v_j$
are not convergent in $H^{s+4}_{\d}$, then from (3.12) we see that
$L_1 L_1^*$ has closed range. But this is a contradiction, since
$L^*_1$ has closed range and, for any operator $A$,
$AA^*$ has closed range if and only if $A$ does.

This completes the proof of Proposition (3.11). \qed
\enddemo

The second part of the uniform surjectivity is the
\proclaim{(3.14) Proposition} Suppose that $f \in H^s_{-\d}$ and
let $v+w$ be the (unique) solution of $L_T (v+w) = f$ in
$H^{s+2}_{-\d} \oplus W$. Then the function $\phi$ in the
nullspace of $L_T$ in $H^{s+2}_\d$ defined by
$\phi = (v+w) - G_T f$ satisfies $\|\phi\| \le C \|f\|_{s,-\d}$.
\endproclaim
\flushpar We have not specified which space the norm of $\phi$ is
to be taken. However, since $\Cal B$ is finite dimensional,
all choices are equivalent. To be definite, we take it as the
$L^2$ norm over a compact set $K = K_1 \cup K_2$ where $K_i \subset M_i
\setminus B_{\a_i}(p_i)$.

\demo{Proof:} Again suppose not, so that for some sequence $T_j$
tending to infinity there is an element $f_j \in H^s_{-\d}$
with $\|f_j\|_{s,-\d} \rightarrow 0$, such that the corresponding
$\phi_j$ satisfies $\|\phi_j \| = 1$. Then we may take a limit
and get a nontrivial function $\phi$ on $M_1 \setminus \{p_1\}$
which is bounded on the cylindrical end and satisfies $L_{1,c} \phi= 0$.
Then, as before, $\psi = u_1^{-1} \phi$ is smooth across $p_1$
and solves $L_1 \psi = 0$. But since $\phi_j = v_j + w_j - G_{T_j}f_j$
and $\|f_j\|$ tends to zero, we see that $\psi \in H^{s+2}_{-\d}
\oplus W_1$, which is a contradiction since we assumed that
$L_1$ has no nullspace here. \qed
\enddemo

\head Proof of Theorem (3.1) \endhead

It is now a relatively simple matter to complete the proof of the
gluing theorem.  In fact, using either (2.15) or (2.17) as appropriate,
the nonlinear step is trivial. Recall that we wish to find a
small function $v \in H^{s+2}_\d$ or pair $v = (\tilde v,w)$ near zero in
$H^{s+2}_{-\d} \oplus W$ such that
$$
N_{g_{c,T}}(u_T + v) = 0.  \tag 3.15
$$
We shall suppose that we are in the case where (2.17) applies,
to be definite.
If the $w$ component were trivial, this would be equivalent
to requiring that $(1+\tilde v)^{4/(n-2)}g_T$ have CPSC, but in
general $w$ cannot be treated as a simple conformal factor,
cf. (2.9). In general, (3.15) is equivalent to the condition
that $(1+ \tilde v)^{4/(n-2)}g^{w}_T$ has CPSC $n(n-1)$, where $g^{w}_T$
is defined in analogous fashion to ${\tilde g}_{\e,a,b}$ in \S 2.

Denote the right side of (3.15) by $N(v)$. Expanding in a Taylor
series shows that
$$
N(v) = f_T + L_T v + Q_T(v),
$$
where $f_T$ is the error term (3.3), $L_{T}$ is the linearization
of $N$ at $u_T$, hence acts by $L_T(v) = L_T(\tilde v + w)$,
and $Q_T$ is a quadratically small remainder term in $v$, which
depends uniformly on $u_T$, hence on $T$.
The inverse $\tilde G_T$ of $\tilde L_T$ is an isomorphism from
$H^s_{-\d}$
to $H^{s+2}_{-\d}\oplus W$, and so we can rewrite (3.15) now as
$$
v = - G_T \left( f_T + Q_T(v) \right), \qquad v \in V. \tag 3.16
$$

It is now standard to solve (3.16), for example by showing that the map
$$
\eta \longrightarrow \Cal T(\eta) = - G_T \left( f_T + Q_T(\eta) \right)
$$
is a contraction on a ball of radius $\sigma$ about $0$ in
$H^{s+2}_{-\d}\oplus W$, hence has a unique fixed point.  In fact,
by taking $T$ large enough we have $\|f_T\| \le C e^{-T/2} \le \sigma/2A$
and $\|Q_T(\eta)\| \le C \sigma^2 \le \sigma/2A$ for $\|\eta\| \le \sigma$,
where $A$ is a uniform bound for the norm of $G_T$. Thus
$\Cal T$ maps the ball of radius $\sigma$ into itself. Furthermore,
$$
\|\Cal T(\eta_1) - \Cal T(\eta_2)\| \le A\|Q_T(\eta_1) - Q_T(\eta_2) \| \le
2A\sigma \|\eta_1 - \eta_2 \|,
$$
hence if $\sigma < 1/4A$, $\Cal T$ is a contraction mapping and there
is a unique $v\in V$ satisfying (3.16) and (3.15).

The only remaining fact to check is that the solution metric $g$ obtained
here is nondegenerate. By the locality of $W$, (2.16) and (2.17) are
immediate. (2.12) can be proved by contradiction. It suffices
to prove that the linearization of $N$ at $g$ is injective
on $H^s_{-\d}$ provided $T$ is sufficiently large. The proof
of this is identical to that of Proposition (3.5).
This completes the proof of Theorem (3.1). \qed

We remark that the solution $v$ of (3.15) and (3.16) is unique
in $H^{s+2}_{-\d}\oplus W$, but it is not the unique solution near zero
in $H^{s+2}_{\d}$.  In fact, since we are working orthogonal to
$\Cal B\cap (H^{s+2}_{-\d}\oplus W)$, where
$\Cal B\equiv \ker (L) \cap H^{s}_{\d}$, we see that the solutions in
$H^{s+2}_{\d}$ are parametrized, as in \cite{MPU}, by the elements of
$\Cal B$.

\specialhead IV. Dipole metrics and applications to the moduli spaces
\endspecialhead

In this section we wish to discuss the ramifications of Theorem (3.1)
in the special case where the component manifolds $(M_i,g_i)$ are
conformally cylinders with Delaunay metrics. As established in \S 2,
these solutions are nondegenerate, and the deficiency
subspace $W$ may be localized to one end of each cylinder.
Furthermore, it is clear that this gluing procedure
can be iterated, because the glued solutions are again nondegenerate.
This leads us to
\proclaim{(4.1) Theorem} Let $(M_i, g_{\e_i})$, $i = 1, 2, \dots$,
be any sequence of Delaunay manifolds. Choose points
$p_i^+,\,p_i^-$ and sufficiently small balls $B_{\a_i^\pm}(p_i^\pm)$
on each $M_i$. Then for each $N \ge 2$, and for gluing parameters
$\eta_j$ sufficiently small, there is a nondegenerate metric
$g^{(N)}(\eta)$ with
CPSC on the iterated connected sum
$$
M^{(N)}(\eta) \equiv M_1 \#_{\eta_1} M_2 \#_{\eta_2} \dots \#_{\eta_{N-1}}M_N,
$$
obtained by gluing the neighbourhood $B_{\a_i^+}(p_i^+)$ to
$B_{\a_{i+1}^-}(p_{i+1}^-)$, which is near to the connected sum metric.
This metric may be obtained inductively, by gluing
$(M_N, g_{\e_N})$ to $(M^{(N-1)}(\eta'),g^{(N-1)}(\eta'))$, where
$\eta' = (\eta_1, \dots, \eta_{N-1})$. Furthermore, if the
sequence $\{\eta_j\}$ tends to zero rapidly enough, then the
Riemannian manifolds $(M^{(N)}(\eta),g^{(N)}(\eta))$ converge on
compact sets to a manifold $M$ with infinitely generated
homology group $H^{n-1}$, and a Riemannian metric
$g$on $M$ of CPSC.
\endproclaim

\demo{Proof}
The only statement that needs proving is the convergence as $N \rightarrow
\infty$. By the Harnack inequality, it suffices to show that
in this iterative process,
the norm of the conformal factor does not blow up or degenerate
on some fixed compact set $K \subset M_1$. But in the gluing
procedure we have shown that we can make the conformal factor
as small as desired on the set $K$ by choosing $T$ sufficiently
large. Thus, when gluing $M_N$ onto the previously constructed
connected sum, we choose $T_N$ so that the conformal factor on
$K$ is no larger than $2^{-N-2}$. The net change of all the
conformal factors is then no larger than $1/2$, hence the
procedure converges. \qed
\enddemo

In \cite{MPU} we studied the moduli space $\ML$ of all CPSC metrics $g$ on
the complement of a fixed singular set $\La$ in $S^n$.
The basic result is that $\ML$ is a real analytic set.
As such, it is endowed with a real analytic stratification
into smooth analytic varieties, but without further
information it is impossible to control the dimensions of
these strata or of the whole moduli space. In particular,
it is conceivable that $\ML$ could be a single point.
However, if one can prove that some $g \in \ML$ is
nondegenerate, then a neighbourhood of any nondegenerate solution $g$
is a real analytic manifold of dimension $\# \La$. This shows
that the connected component of $\ML$ containing $g$ is of this dimension,
and this nondegenerate solution is in the principal stratum.
Unfortunately, it was not at all clear that any of those
solutions previously known to exist, i.e. those constructed by Schoen,
are nondegenerate. The original motivation for this paper
was to construct nondegenerate solutions, at least for
certain singular sets $\La$, and this has now been accomplished.
Thus we have proved

\proclaim{(4.2) Corollary} Let $(M^{(N)}, g^{(N)})$ be one of the solutions
obtained in Theorem (4.1) by gluing together $N$ cylinders with their
Delaunay metrics. This manifold is conformally flat and may be
uniformized as a domain in the sphere $\Omega = S^n \setminus \La$, for
some set $\La$ with $2N$ elements. The resulting conformally
flat metric $g = g^{(N)}$ is a nondegenerate element of the moduli space
$\ML$.  The component of $\ML$ containing $g$ is a $2N$-dimensional
real analytic set and $g$ is in the principal, top-dimensional,
stratum.
\endproclaim

The singular sets $\La$ we obtain by this gluing procedure
are very special: they contain $N$ pairs of closely spaced
points, with the distance between each pair bounded from below.
We call these dipole configurations. In the next subsection
we indicate how we may use the nondegeneracy of dipole metrics
to deduce nondegeneracy of elements in other moduli spaces $\ML$
where $\La$ is not necessarily a dipole configuration.

We end this subsection by observing that by the construction
we can prescribe the geometry of half of the ends of a dipole solution.
This follows immediately by observing that the elements in the parameter
space $W$ are supported on only half (i.e. $\{(t_j,\theta_j):t_j>-1\}$)
of each  cylinder.

\proclaim{(4.3) Corollary} Let $g$ be a dipole metric constructed
by gluing together the $N$ Delaunay solutions $D_{\e_1}, \dots,
D_{\e_N}$. Let the singular set $\La \subset S^n$ be written
as a set of pairs $\{(x_1,y_1), \dots, (x_N,y_N)\}$. Then the asymptotic
Delaunay parameter for $g$ at $x_j$ is exactly $\e_j$ and the Delaunay
parameter at $y_j$ is very close to $\e_j$. In particular, since we
can prescribe the Delaunay parameters at $x_j$, by choosing $\e_j=\bar u$
for $j=1,\ldots,N$,  there exist dipole solutions
with $N$ asymptotically cylindrical ends.
\endproclaim

This shows that the dipole solutions are geometrically very different
than the solutions constructed by Schoen, where the Delaunay
parameters at each singular point are very close to zero.

\head The unmarked moduli space ${\Bbb M}_k$ \endhead

As a final topic, and to further elucidate the impact of the dipole
solutions on the moduli space theory, we shall consider
the unmarked moduli space
$$
\Mk = \{(g,\La):\La\in \Ck \text{ and }g\in\ML\}
$$
of conformally flat CPSC metrics on the complement of any
$k$ points in the sphere.  Here $\Ck$ is the configuration space of
$k$ distinct points in $S^n$ and $\ML$ is the moduli space of solutions
on $S^n \setminus \La$ studied in \cite{MPU}. There is a natural map
$$
\pi: \Mk \longrightarrow \Ck.
$$

In this subsection we prove that $\Mk$ is a real analytic set.
We will only sketch the proof since it is very close to the
proof of the analogous result in \cite{KMP}; note that
we use this rather than the proof in \cite{MPU} because
it is much more direct and simple to modify.
A consequence of this result is that $\Mk$ admits a stratification
into smooth real analytic manifolds. The dimensions of these
strata, particularly the maximal dimension, are difficult to
obtain without more information. However, there is a
nondegeneracy condition which is less stringent than the
one we have used before, such that if $(g,\La) \in \Mk$ is
nondegenerate in this new sense, then the stratum containing
this point is maximal (in that connected component of $\Mk$) and has
the expected dimension $k(n+1)$.  We call this new nondegeneracy
condition unmarked nondegeneracy, and the former one simply nondegeneracy,
as before.
We shall see that nondegenerate points in $\ML$ are a fortiori nondegenerate
in the unmarked moduli space $\Mk$. In particular, since the dipole metrics
are nondegenerate in $\ML$, any component of $\Mk$, $k = 2N$, containing a
dipole metric has top stratum of the predicted dimension.  Furthermore
all other points in this top dimensional stratum, including ones where the
singular points are no longer in the dipole configurations, are then
nondegenerate in $\Mk$. We shall show, finally, that most of these are
nondegenerate in their respective moduli spaces $\ML$.

Before we can give the precise statements and proofs, we digress briefly
to discuss some additional facts about Jacobi fields for Delaunay metrics
which we need. We also quote several results and the (adaptations of the)
relevant lemmas from \cite{MPU} which are required.

We have already used the two temperate Jacobi fields $\phi_0^\pm$ on
$D_\e$, one periodic and the other linearly growing, which correspond
to infinitesimal translations and changes of Delaunay parameter.
There is another set of Jacobi fields, $\phi_j^\pm$, $j = 1, \dots, n$,
on $D_\e$ which can be written explicitly. Using cylindrical coordinates,
and writing the Delaunay metric $g_\e = u_\e^{4/(n-2)}(dt^2 + d\th^2)$,
then
$$
\phi_j^{\pm} = e^{\pm t}\left( \frac{n-2}{2}u_\e \pm
\frac{du_\e}{dt}\right)\psi_{j}(\theta), \tag 4.4
$$
where $\psi_{j}(\theta)$ are the eigenfunctions for
$\Lap_{\theta}$  with eigenvalue $n-1$ on
$S^{n-1}$, are also solutions of $L_\e \phi = 0$.  Like the $\phi_0^\pm$,
but unlike all the other Jacobi fields for $g_\e$, these arise
as derivatives of explicit one-parameter families of solutions.
In fact, regarding Delaunay solutions as metrics on $S^n$ with two
singular points, we can pull these metrics back by any conformal
transformation of $S^n$ to obtain new solution metrics, also in
$\Bbb M_2$. Families
of solutions are obtained by pulling back by families of
conformal transformations. If we differentiate the family obtained
from the one-parameter family of conformal transformations fixing both
singular points, we obtain $\phi_0^+$. If we differentiate, on
the other hand, the families of parabolic conformal transformations fixing
one or the other of the singular points of the original $g_\e$,
we obtain these other $\phi_j^\pm$. These parabolic transformations
correspond to the translations on the images, $\Bbb R^n$, of the stereographic
projection from $S^n$ with one of the singular points removed.  The parabolic
transformations moving one singular point lead to Jacobi fields blowing up
exponentially at that singular point and decaying exponentially
at the other. Their exact decay rate is given by (4.4).
Combining this with the earlier discussion of the Fredholm
properties of $L_\e$ on weighted Sobolev spaces, we see that $L_\e$ is not
Fredholm on the space $H^s_{\d{_1}}$, with exponential weight
$\d_{1} = P_\e$.  Actually, as proved in \cite{MPU}, $L_\e$ is Fredholm
on all other weighted spaces $H^s_\d$ for $0 < \delta < \d_1+\eta$,
$\delta \ne \d_1$, where $\eta$ is some small positive number depending
on $\e$.

If $g$ is any CPSC metric on $S^n \setminus \La$, with $\La$
finite, then in a neighbourhood of each $x_j \in \La$, $g$ is
asymptotically equivalent to a Delaunay metric. Approximate Jacobi fields
may be constructed by transplanting cutoffs of the Jacobi fields
$\phi_j^\pm$, $j = 0, \dots, n$, from the appropriate $D_\e$ to
a neighbourhood of each $x_j$. This, together with the approximate
Jacobi fields corresponding to $\phi_0^\pm$ yields a $2k(n+1)$ dimensional
family of approximate Jacobi fields, which we denote $\Cal W$; it
substitutes for the deficiency space $W$ in the analysis of the
unmarked moduli space. Note that $\Cal W \subset H^s_{\d_1+\d}$
for any $\d > 0$. Since there are good parametrices
for $L_g$ (constructed in \cite{MPU}), and since the elements
of $\Cal W$, in particular those of $W$, correspond to
one-parameter families of solutions, the marked nondegeneracy
conditions (2.13) and (2.15-17) reduce simply to the hypothesis
that $L_g$ is injective on $H^s_{-\d_1-\delta}$ for any $\delta > 0$.
By contrast, we make the

\proclaim{(4.5) Definition} The solution $(g,\La) \in \Mk$ is called
unmarked nondegenerate if the linearized operator $L_g$
has no nullspace in $H^s_{-\d_1-\d}$ for any $\d > 0$.
\endproclaim

By duality, $g$ is unmarked nondegenerate if and only if
$$
L_g: H^{s+2}_{\d_1+\d} \longrightarrow H^s_{\d_1+\d}
$$
is surjective for $\d>0$ sufficiently small. Using the parametrix
construction and the proof of the linear decomposition lemma (4.18)
from \cite{MPU} we obtain

\proclaim{(4.6) Lemma} Suppose that $(g,\La) \in \Mk$ is unmarked
nondegenerate.
Then there exists a bounded map $G$ from $H^s_{\d_1+\d}$ to
$H^{s+2}_{\d_1+\d}$ such that $L_g G = I$. If we restrict the
domain of $G$ to $H^s_{-\d_1-\d}$, then its range is
$H^{s+2}_{-\d_1-\d} \oplus \Cal W$. In particular
$$
L_g: H^{s+2}_{-\d_1-\d} \oplus \Cal W \longrightarrow H^s_{-\d_1-\d} \tag 4.7
$$
is surjective.
\endproclaim

The dimension of the space $\Cal W$ is $2k(n+1)$ and these form the
possible parameters for deformations of $(g,\La) \in \Mk$.
Which deformations actually occur is a difficult question, however
by a relative index theorem  as in \cite{MPU} we have

\proclaim{(4.8) Lemma} The kernel of the map (4.7) has dimension k(n+1).
\endproclaim

We also require an analogue of (2.17), that the elements of
$\Cal W$ (which are sufficiently small in norm)
can be `exponentiated' to one-parameter families
of metrics which are locally of CPSC near the singular points.
Thus let $w \in \Cal W$; we can write it as a sum
$$
w = \sum_{\ell = 1}^k \sum_{j = 0}^n a_{j,\pm}^{(\ell)}
\chi^{(\ell)} \phi_j^\pm
$$
near each singular point $x_\ell$. Here $\chi^{(\ell)}$ are
fixed cutoff functions equaling one in a small ball around
$x_\ell$ and vanishing outside a slightly larger ball.
We identify $w$ with the collection of coefficients
$\bold a = \{a_{j,\pm}^{(\ell)}\} \in \RR^{2k(n+1)}$.
Then by a procedure very similar to the one in \S 2, we
define a metric $g_{\bold a}$ by altering $g$ in
the neighbourhoods $B_\sigma(x_\ell)$ according to the
amounts specified by the coefficients $\bold a$. The first
step is to identify the Delaunay solution $D_{\e_\ell}$ to
which $g$ is asymptotic in $B_\sigma(x_\ell)$. In this
ball we can write $g = (1 + v)^{4/(n-2)}g_{\e}$ where
$v = O(r)$. We first conformally deform $v$ to equal
zero in $B_{\sigma/2}(x_\ell)$. The new metric $\tilde g$
is exactly Delaunay in these smaller balls, and is unchanged
outside the larger balls. Now we can alter $\tilde g$ in
this smaller ball in a manner specified by the parameters
$\bold a\in \RR^{2k(n+1)}$, the new metric we obtain is denoted
$g_{\bold a}$. We demonstrated earlier how to do this for
the coefficients $a_{0,\pm}=\{a,b\}$ for the Delaunay solutions,
and since this change was supported in the half-cylinder that discussion
obviously localizes to these half-Delaunay solutions.
The analogous procedure for the other coefficients is  defined
exactly as before. Note that the singular points $x_\ell$ will be
moved in this procedure if any one of the $a_{j,+}^{(\ell)}$
is nonzero.

Now we may define the nonlinear operator
$$
N: H^{s+2}_{-\d_1-\d} \oplus \Cal W \longrightarrow H^s_{-\d_1-\d}
$$
for small elements $w \in \Cal W$ at the metric $g \in \Mk$ by
$$
N(v,w) = \Lap_{g_{\bold a}}(1+v) - \frac{n-2}{4(n-1)}R(g_{\bold a}) (1+v)
+ \frac{n(n-2)}{4} (1+v)^{\frac{n+2}{n-2}},
$$
where $w \in \Cal W$ determines the coefficients $\bold a$.
This is clearly a real analytic mapping. A neighbourhood
$\Cal V$ of $g$ in $\Mk$ is given as the zero set
of this map $N$ in $H^{s+2}_{-\d_1-\d} \oplus \Cal W$. Its linearization
coincides with the linearization $L_g$ acting on $H^{s+2}_{-\d_1-\d} \oplus
\Cal W$. The analytic implicit function theorem and the relative index
theorem of [MPU] now give

\proclaim{(4.9) Proposition} If $g \in \Mk$ is unmarked nondegenerate,
so that
$$
L_g: H^{s+2}_{-\d_1-\d} \oplus \Cal W \longrightarrow H^s_{-\d_1-\d}
$$
is surjective, then there is a neighbourhood $\Cal V$ of $g \in \Mk$
which is a real analytic manifold of dimension $k(n+1)$.
\endproclaim

When $g$ is unmarked degenerate a somewhat weaker statement is
true:

\proclaim{(4.10) Proposition} If $g \in \Mk$ is unmarked degenerate,
then there is a neighbourhood $\Cal X$ of $(0,0) \in H^{s+2}_{-\d_1-\d}
\oplus \Cal W$ and a real analytic diffeomorphism $\Psi$ from
this neighbourhood to itself such that $\Psi(N^{-1}(0)\cap \Cal X)$ lies
in a neighbourhood of zero $\Cal Y$ in finite dimensional subspace
$P$ and coincides with the zero set of a real analytic
mapping from $\Cal Y$ to $\CC^N$, where $N$ is the dimension
of the kernel of $L_g$ on $H^{s}_{-\d_1-\d}$.
\endproclaim

The proof of this final proposition is identical to the one
in \cite{KMP}, and thus we give only the briefest sketch and refer
there for the details.
We use what is often called the Liapunov-Schmidt or Kuranishi method.
Let $K$ denote the nullspace of $L$ on $H^s_{-\d_1-\d}$, which
is nontrivial since $g$ is unmarked degenerate. Then
$$
\Cal L: H^{s+2}_{-\d_1-\d} \oplus \Cal W \oplus K \longrightarrow
H^s_{-\d_1-\d},
$$
defined by $\Cal L(v,w,k) = L_g(v+w) + k$ is surjective,
by construction. We can similarly define
$\Cal N(v,w,k) = N(v,w) + k$, so that $\Cal L$ is the (surjective)
linearization of $\Cal N$. The zero set of $N$ is identified with the
set $\{(v,w,k): \Cal N(v,w,k) = k\}$. Using the implicit function
theorem for $\Cal N$, we can show that this set is real analytic.

It is clear that the projection $\pi: \Mk \longrightarrow \Ck$
is real analytic.

We now combine these last two results with the results of the
last subsection.

\proclaim{(4.11) Corollary} Suppose that $k = 2N$ and that
some component of $\Mk$ contains an element $(g,\La)$ where
$\La$ is a dipole configuration and $g$ a dipole metric.
Then $(g,\La)$ is in the principal stratum of that component,
and this principal stratum is nonsingular and has dimension $2N(n+1)$.
Let $\tilde \pi$ denote the restriction of the projection $\pi$
to this component of $\Mk$ and let $\Cal C$ be the image
of $\tilde \pi$ in $\Ck$. Then $\Cal C$ is a subanalytic
set, with principal stratum of dimension $kn=2Nn$, and
the preimage $\Cal M_{\La'}$ of any configuration $\La'$ in this principal
stratum is a moduli space with a nonsingular $k=2N$-dimensional principal
stratum composed of nondegenerate elements.
\endproclaim

As a final remark we note that the singular strata for these moduli
spaces may in fact be trivial. In \cite{MPU} we show that they
are trivial for generic conformal classes. In fact, we do not know
whether degenerate solutions exist in any case. Either a
construction of degenerate solutions, or a geometric criterion
for their presence or absence, would be quite interesting.

\vfill\eject

\enddocument


\Refs

\widestnumber\key{KKMS}

\ref
\key D
\by C. Delaunay
\paper Sur la surface de revolution dont la courbure moyenne est constant
\jour J. Math. Pure Appl.
\vol 6
\yr 1841
\pages 309-320
\endref

\ref
\key E
\by J. F. Escobar
\paper Conformal deformation of a Riemannian metric to a constant
scalar curvature metric with constant mean curvature boundary
\paperinfo Preprint
\endref


\ref
\key F1
\by R. H. Fowler
\paper The form near infinity of real continuous solutions of a certain
differential equation of the second order
\jour Quarterly J. of Pure and  Appl. Math.
\vol 45
\yr 1914
\pages 289-349
\endref


\ref
\key F2
\by R. H. Fowler
\paper Further studies of Emden's and similar differential equations
\jour Quarterly J. of Math., Oxford series
\vol 2
\yr 1931
\pages 259-287
\endref


\ref
\key J
\by D.D. Joyce
\paper Hypercomplex and Quaternionic Manifolds and
Scalar Curvature on Connected Sums
\paperinfo D.Phil. thesis, Oxford University
\yr 1992
\endref

\ref
\key K
\by N. Kapouleas
\paper Complete constant mean curvature surfaces in Euclidean three-space
\jour Ann. of Math.
\vol 131
\yr 1990
\pages 239-330
\endref


\ref
\key KMP
\manyby  R. Kusner, R. Mazzeo and D. Pollack
\paper The moduli space of complete embedded constant mean curvature surfaces
\paperinfo To appear Geom. and Funct. Analysis
\endref

\ref
\key MPU
\manyby R. Mazzeo, D. Pollack and K. Uhlenbeck
\paper Moduli spaces of singular Yamabe metrics
\paperinfo To appear, Jour.  Amer. Math. Soc
\endref

\ref
\key M1
\by R. Mazzeo
\paper Regularity for the singular Yamabe equation
\jour Indian Univ. Math. J.
\vol 40
\yr 1991
\pages 1277-1299
\endref

\ref
\key M2
\by ---
\paper Elliptic theory of differential edge operators I
\jour Comm. in P.D.E.
\vol 16
\yr 1991
\pages 1615-1664
\endref

\ref
\key MPa1
\by --- and F. Pacard
\paper A new construction of singular solutions for a semilinear
elliptic equation
\paperinfo To appear, J. Differential Geometry
\endref

\ref
\key MPa2
\by --- and ---
\paper Singular Yamabe metrics with isolated singularities
\paperinfo In preparation
\endref

\ref
\key MS
\manyby --- and N. Smale
\paper Conformally flat metrics of constant positive scalar curvature
       on subdomains of the sphere
\jour J. Differential Geometry
\vol 34
\yr 1991
\pages 581-621
\endref


\ref
\key P1
\by D. Pollack
\paper Nonuniqueness and high energy solutions for a conformally invariant
       scalar equation
\jour Comm. Anal. and Geom.
\vol 1
\yr 1993
\pages 347-414
\endref

\ref
\key P2
\by ---
\paper Compactness results for complete metrics of constant positive scalar
       curvature on subdomains of $S^{n}$
\jour Indiana Univ. Math. J.
\vol 42
\yr 1993
\pages 1441-1456
\endref

\ref
\key S1
\by R. Schoen
\paper The existence of weak solutions with prescribed singular
       behavior for a conformally invariant scalar equation
\jour Comm. Pure and Appl. Math.
\vol XLI
\yr 1988
\pages 317-392
\endref

\ref
\key S2
\by ---
\paper Variational theory for the total scalar
       curvature functional for Riemannian metrics and related topics
\inbook Topics in Calculus of Variations; Lecture Notes in Mathematics
\# 1365
\publ  Springer-Verlag
\eds M. Giaquinta
\yr 1987
\pages 120-154
\endref


\ref
\key SY
\manyby --- and S.T. Yau
\paper Conformally flat manifolds, Kleinian groups and scalar curvature
\jour Invent. Math.
\vol 92
\yr 1988
\pages 47-71
\endref

\ref
\key T
\by C. Taubes
\paper Self-dual Yang-Mills connections on non-self-dual 4-manifolds
\jour  J. Differential Geometry
\vol 17
\yr 1982
\pages 139-170
\endref


\endRefs
\bye